\documentclass{article}

\usepackage{amsfonts,amsmath,amssymb}

\newcommand{\E}{\mathsf{E}}
\newtheorem{lemma}{Lemma}

\newcommand{\Rc}{\mathcal{R}}
\newtheorem{theorem}{Theorem}

\newtheorem{remark}{Remark}
\newcommand{\Pp}{\mathsf{P}}
\newcommand{\R}{\mathbb{R}}

\newcommand{\N}{\mathbb{N}}

\newcommand{\Bc}{\mathcal{B}}
\newcommand{\Dc}{\mathcal{D}}
\newcommand{\C}{\mathcal{C}}
\newcommand{\ONE}{{\bf 1}}

\newcommand{\kap}{{\varkappa}}
\newcommand{\bpf}[1][Proof]{{\noindent {\sc #1: }}}
\newcommand{\epf}{{{\hspace{4 ex} $\Box$ \smallskip}}}
\makeatletter
\def\imod#1{\allowbreak\mkern1mu({\operator@font mod}\,\,#1)}
\makeatother

\author{Yuri Bakhtin\thanks{School of Mathematics, Georgia Tech, Atlanta GA, 30332-0160;  email:bakhtin@math.gatech.edu, 404-894-9235 (office phone), 404-894-4409(fax)}}

\title{Self-Similar Markov Processes on Cantor Set}
\begin{document}
\maketitle
\begin{abstract}
 We define analogues of Brownian motion on the triadic Cantor set by introducing a few natural requirements on the Markov
semigroup. We give
a detailed description of these symmetric self-similar processes and study their properties such as mixing and moment asymptotics.
\end{abstract}

\section{Introduction}
In a recent paper \cite{Pearson-Bellisard}, the authors use noncommutative geometry to describe an analogue of the Riemannian structure for 
ultrametric Cantor sets. This eventually leads them to a definition of a certain family of operators that play the role of the Laplace--Beltrami operators for Cantor sets. It is natural to treat these operators as infinitesimal generators of Markov semigroups on Cantor sets and call the associated Cantor-set-valued stochastic processes to be the analogues of the Brownian motion.

The goal of this note is to provide an alternative definition of Brownian motion on the classical triadic Cantor set. We use the axiomatic method
and describe several natural requirements, most important of which are isometry invariance and scale invariance, that should hold for a reasonable analogue of Brownian motion on the Cantor set. Then we give a complete description of Markov processes satisfying these requirements. We call these processes symmetric self-similar (SSS).

The parametrization of SSS processes involves two degrees of freedom. One of these is responsible just for uniform time changes, so that effectively this family of SSS processes is parametrized by one parameter of scaling, or self-similarity.
Since we are basing our approach on scaling properties, the SSS processes are, in fact, analogous to
symmetric stable L\'evy processes in~$\R$. It is remarkable that the class of SSS processes exactly coincides with the class of processes on the triadic Cantor set obtained in~\cite{Pearson-Bellisard} as a result of a completely different
approach.

Our approach is somewhat similar to Schramm's celebrated characterization of SLE via conformal invariance and 
Markov property, see e.g.~\cite[Chapter 6]{Lawler:MR2129588}.

It is worth mentioning that random walks on self-similar fractals have been studied in the literature, see e.g.~\cite{Barlow:MR1668107}. However, to the best of our knowledge, our results for a {\it disconnected ultrametric} Cantor set are new. We also stress that our approach can be easily implemented for other Cantor sets with rich self-semilarity properties. However, it is not clear how it can be used for a general Cantor set.

The paper is organized as follows. In Section~\ref{sec:ultrametric} we introduce the setting by a description of the geometry of the Cantor set. In Section~\ref{sec:main_characterization_theorem}, we give our main result, a complete characterization of SSS Markov processes on the Cantor set, and briefly discuss the relationship with the processes studied in~\cite{Pearson-Bellisard}. In Section \ref{sec:explicit_construction} we describe an explicit construction of SSS processes and study their jump
statistics. The role of the Laplacian on the Cantor set is played by the generators of SSS processes, and in 
Section~\ref{sec:spectral_structure} we describe some of their spectral properties. In Section~\ref{sec:ergodicity},
we prove that the Bernoulli measure on~$C$ is a unique invariant measure for any SSS process. We also prove that it is exponentially attracting. In Section~\ref{sec:moments}, we study the asymptotics of displacement moments
of SSS processes for small transition times.

{\bf Acknowledgements.} The author is grateful to Jean Bellissard and John Pearson for introducing him to their work on the geometry of Cantor sets, for several stimulating discussions that led to writing this paper, and for some corrections.
He is also grateful to NSF for partial support of this research via CAREER award DMS-0742424.

\section{The ultrametric structure on the triadic Cantor set}\label{sec:ultrametric}
A Cantor set  is a topological space that is non-empty, compact, perfect,
totally disconnected and metrizable. In this paper, we study a classical example, the triadic Cantor set 
\[
C=\left\{x:\ x=\sum_{k=1}^\infty \frac{2x_k}{3^k},\  x_k=0,1,\ \mbox{\rm for all}\ k\in\N\right\}.
\]
In many situations it is natural to identify $x\in C$ with the sequence $(x_k)_{k\in\N}$ which in turn may be identified with
an infinite simple path on the infinite rooted binary tree.

For every $x,y\in C$ we define 
\[d(x,y)=3^{-c(x,y)},\]
 where 
\[
c(x,y)=\min\{k\in\N: x_k\ne y_k\}.
\]
(We agree that $c(x,x)=\infty$ and $d(x,x)=0$ for any $x\in C$.)
 It is easy to see that 
$d$ is equivalent to the Euclidean metric:
\[
\frac13|x-y|\le d(x,y)\le |x-y|,\quad x,y\in C.
\]
It is also easy to see that 
$d$
is an ultrametric, i.e.\ it satisfies the strong triangle
inequality: 
\[
d(x, y)\le\max\{d(x, z), d(z, y)\},\quad x,y,z\in C.
\]

The metric space $(C,d)$ has rich structure that involves self-similarity properties and a rich group of isometries. 
Let us denote
\[
 \pi_nx=(x_1,\ldots,x_n),\quad x\in C, n\in\N.
\]

Then, for every $n\in\N$,  the set $C$ may be decomposed into a disjoint family of $2^n$ sets:
\[
[v]=\{x\in C:\ \pi_n=v\},\quad v\in L_n,
\] 
where $L_n=\{0,1\}^n$.

Each of these sets is similar to $C$ with scaling coefficient $3^n$. One of the similarity maps is given by
\[
(v_1,v_2,\ldots,v_n, x_{n+1},x_{n+2},\ldots)\mapsto  (x_{n+1},x_{n+2},\ldots).
\] 
 
In particular, these sets are isometric to each other and have diameter $3^{-n-1}$. This leads to the following complete description of all isometries. Let $g$ be an isometry.
Then for each $n\in\N$, $g$ permutes sets  $[v],\ v\in L_n$, i.e.\ for every $v\in L_n$, there is $v'\in L_n$
such that
\[
g([v])=[v'], 
\]
which we shall abbreviate as
\[
 g_n(v)=v'. 
\]

Obviously,
\begin{equation}
\label{eq:isometry_consistency}
g_n(v)=v' \Rightarrow  g_{n-1}(\pi_{n-1} v)=\pi_{n-1}v',
\end{equation}
where we used the notation
\[
 \pi_m(v_1,\ldots,v_n)=(v_1,\ldots, v_m),\quad m\le n,\ v\in L_n,
\]
and one can easily show that a sequence of bijections $g_n:L_n\to L_n$ satisfying the consistency
condition~\eqref{eq:isometry_consistency} generates an isometry.

\section{Markov processes, the main characterization theorem}\label{sec:main_characterization_theorem}
Our aim is to obtain an analogue of the Brownian Motion (or symmetric stable Levy process) on $(C,d)$. Since $C$ is completely disconnected, there is no hope that a nontrivial Markov process will have continuous trajectories. So we shall relax the continuity requirement, and consider stochastically continuous Markov processes
with c\`adl\`ag (right-continuous and with left limits) trajectories.
 Let us recall
that a homogeneous Markov processes with transition probability function 
\[
P(t,x,B),\quad t\in[0,\infty),\ x\in C,\ B\in\Bc(C,d), 
\]
where $\Bc(C,d)$ is the Borel $\sigma$-algebra on $(C,d)$, is called stochastically continuous if
for any open set $U$ and any point $x\in U$,
\[
\lim_{t\downarrow 0}P(t,x,U)=1.
\]

Transition probability function $P(\cdot,\cdot,\cdot)$ is called
Feller if the space $\C(C)$ of continuous real-valued functions on $C$
is invariant under the semigroup $(S^t)_{t\ge0}$ generated by $P(\cdot,\cdot,\cdot)$ and defined via:
\begin{equation}
 (S^t f)(x)=\int_C f(y)P(t,x,dy).
\label{eq:semigroup}
\end{equation}
(We refer to Chapter I of~\cite{Liggett:MR776231} for a concise exposition of the necessary background on Markov semigroups.)

Let $(P_x)_{x\in C}$ be a Markov family  compatible with $P(\cdot,\cdot,\cdot)$, i.e.\ 
for each~$x$, $P_x$ is a measure on the space of c\`adl\`ag $C$-valued trajectories on $[0,\infty)$ such that for the
canonical process $X$,
$P_x\{X(0)=x\}=1$ and, under $P_x$, $X$ is a Markov process with transition probabilitiy $P(\cdot,\cdot,\cdot)$. 

We are going to impose restrictions on $P(\cdot,\cdot,\cdot)$ one by one and show that these restrictions together lead
to a concise analytic description of the generator of $(S^t)_{t\ge0}$ allowing for a parametrization of the set of allowed processes by two parameters. Then, for any choice of these two parameters,  we construct a unique process with required properties.

We say that $P(\cdot,\cdot,\cdot)$ is invariant under isometries if for any isometry $g$,
\[
 P(t,g(x),g(B))=P(t,x,B),\quad t\in[0,\infty),\ x\in C,\ B\in\Bc(C,d).
\]

\begin{lemma}
Let $(X_t)_{t\ge0}$ be a homogeneous Markov process with transition function $P(\cdot,\cdot,\cdot)$ invariant under isometries.
Then $(\pi_nX_t)_{t\ge0}$ is also a homogeneous Markov process. 
\end{lemma}

\bpf
We notice that  the isometry invariance straightforwardly implies
\begin{multline*}
\Pp\{\pi_nX_{0}=v_0, \pi_n X_{t_1}=v_1,\ldots, \pi_n X_{t_m}=v_m\}\\
=\int_{[v_0]}\nu(dx_0)\int_{[v_1]}P(t_1-t_0,x_0,dx_1)\int_{[v_2]}P(t_2-t_1,x_1,dx_2)\ldots\int_{[v_m]}P(t_m-t_{m-1},x_{m-1},dx_m)\\
=\nu([v_0])P_n(t_1-t_0,v_0,v_1)\ldots P_n(t_m-t_{m-1},v_{m-1},v_m),
\end{multline*} 
where $\nu$ is the distribution of $X_0$, and the quantities
\[
P_n(s,u,v)=P(s,x,[v]),\quad \pi_nx =u,
\]
are well-defined due to the isometry invariance.

Since $\pi_n(X_t)$ takes finitely many values, we conclude that $\pi_nX$ is a Markov chain with transition matrix
$P_n(s,v,u)$. \epf

The process $\pi_n X$ inherits the stochastic continuity from $X$. Since it takes finitely many values, the
transition rates
\[
q(u,v)=\lim_{t\downarrow 0}\frac{P_n(t,u,v)-\delta_{uv}}{t},
\]
are well-defined, where $\delta_{uv}$ is the Kronecker symbol, see e.g.\ \cite{Chung:MR0217872}[Theorem~5, Section~II.2].
If $v,u\in L_n$ for some $n$, then for any $x\in u$, we set $q(x,v)=q(u,v)$.

Let us introduce 
\[
d(u,v)=d([u],[v])=\inf\{ d(x,y): x\in[u], y\in[v]\}.
\]

Due to the isometry invariance, $q(u,v_1)=q(u,v_2)$ if $d(u,v_1)=d(u,v_2)$, so that 
$q_n(u,v)$ is actually a function of $d(u,v)$, and we can	 write
\[
q(u,v)=q(d(u,v)).
\]
If $n\le m$ and $v\in L_n$ then
\[
q(x,v)=\sum_{\substack{w\in L_m\\\pi_n w=v}}q(x,w),\quad x\in C.
\]
If, moreover, $x\notin [v]$, then considering the family of isometries that leave $x$ fixed and permute points $w$ in the above summation, we can conclude that all the terms $q(x,w)$ coincide. Therefore, in this case,
\begin{equation}
\label{eq:rates_decomposed_as_sums}
 q(x,v)= 2^{m-n} q(x,w),\quad w\in L_m, \pi_n w=v.
\end{equation}

For any $x\in C$ and $n\in\N$, we define $q_n=q(x,v_n(x))$, where
\begin{equation}
\label{eq:v_k}
v_n(x)=(x_1,\ldots,x_{n-1}, 1-x_n).
\end{equation}
 For any $n$,  $q_n$ does not depend on $x$ due to the isometry invariance. Using~\eqref{eq:rates_decomposed_as_sums}, we see that all rates $q(x,v)$ can be expressed
in terms of $q_n$. Namely, for $v\in L_n$ and $x\notin [v]$,
\begin{equation}
\label{eq:jump_rates_uniformly_spread}
q(x,v)=2^{-(n-c(x,v))}q_{c(x,v)},
\end{equation}
where $c(x,v)=\min\{k: x_k\ne v_k\}$. We conclude that the distribution of $X$ under $P_x$ is completely determined
by the family of jump rates $(q_n)_{n\in\N}$.

Let us recall that 
the infinitesimal generator for the semigroup $(S^t)_{t\ge0}$ defined in \eqref{eq:semigroup} is given by
\begin{equation}
Af = \lim_{t\downarrow 0}\frac{S^tf-f}{t}
\label{eq:denerator}
\end{equation}
for $f\in\Dc(A)$, where $\Dc(A)$ is the domain of $A$, i.e., set of all functions $f$ such that the r.h.s.\ of 
\eqref{eq:denerator} is well-defined as a uniform limit.
If $f$ is a cylindric function, i.e., $f(x)=h(\pi_n x)$ for some function $h:L_n\to\R$, then
$f\in\Dc(A)$ and
\begin{equation}
Af (x)=\sum_{v\in L_n}q(x,v) (h(v)-h(\pi_nx)). 
\label{eq:generator1}
\end{equation}

The next property we would like to require is self-similarity, or scale invariance. Let us recall that for any $n$ and any $v\in L_n$, $[v]$ is similar to $C$. 
So we shall require that the distribution of the Markov process $X$ confined to~$[v]$ coincides with the appropriately scaled distribution of the unrestricted process on the entire $C$.

To make this precise, we need to introduce confinements of Markov processes.
For any $T>0$, $n\in\N$, $v\in L_n$, $x\in [v]$ and consider the conditional measure
\[
P^{v,T}_x\left\{X_{[0,T]}\in\cdot\ \right\}=P_x\left\{X_{[0,T]}\in\cdot\ |\ X_t\in [v], t\le T\right\},
\]
where $X_{[0,T]}$ denotes the restriction of the canonical process $X$ onto the time interval $[0,T]$.
\begin{lemma}
\label{lm:restriction_is_markov}  For any $T>0$, $n\in\N$, $v\in L_n$, $x\in [v]$, the canonical process $X$ is Markov
under measure $P^{v,T}_x$ on $[0,T]$. Moreover, these distributions are consistent for different values of $T$. 
\end{lemma}
\bpf 
Choose any $n'>n$, $l\ge 1$, and  $u_1,\ldots,u_l\in L_{n'}$ such that $[u_k]\subset[v]$ for all $k$,
and $0<t_1<\ldots<t_l=T$. Then
\begin{equation}
\label{eq:checking_Markov_property_1}
P_x\bigl\{X_{t_1}\in[u_1],\ldots, X_{t_l}\in[u_l]\ |\ X_t\in[v], t\in[0,T]\bigr\}=\frac{b}{a}.
\end{equation}
Here
\begin{align*}
 a=P_x\{X_t\in [v], t\le T\}&=\lim_{m\to\infty} P_x\left\{X_t\in [v], t=\frac{T}{m},\frac{2T}{m},\ldots,\frac{mT}{m}\right\}
 \\&=\lim_{m\to\infty} P\left(\frac{T}{m},v,v\right)^m
 \\&=\exp\{q(v,v) T\},
\end{align*}
(we used the right-continuity of trajectories for the first identity,   the isometry
invariance plus for the second, and the definition of $q_n$ for the third one), and
\begin{align*}
b&=P_x\{X_{t_1}\in[u_1],\ldots, X_{t_l}\in[u_l]; \ X_t\in[v], t\in[0,T]\}\\
&=\lim_{m\to\infty}\sum  P\left(\frac{t_1}{m},[\pi_n x],[w_1^1]\right)P\left(\frac{t_1}{m},[w_1^1],[w_2^1]\right)\ldots P\left(\frac{t_1}{m},[w_{m-1}^1],[u_1]\right)\\
&\quad\times P\left(\frac{t_2-t_1}{m},[u_1],[w_1^2]\right)P\left(\frac{t_2-t_1}{m},[w_1^2],[w_2^2]\right)\ldots P\left(\frac{t_2-t_1}{m},[w_{m-1}^2],[u_2]\right)\\
&\quad\ldots\\
&\quad\times P\left(\frac{t_l-t_{l-1}}{m},[u_{l-1}],[w_1^l]\right)P\left(\frac{t_l-t_{l-1}}{m},[w_1^l],[w_2^l]\right)\ldots P\left(\frac{t_l-t_{l-1}}{m},[w_{m-1}^l],[u_l]\right)\\
&= \exp\{t_1 Q\}_{\pi_nx,u_1}\exp\{(t_2-t_1) Q\}_{u_1,u_2}\ldots\exp\{(t_l-t_{l-1}) Q\}_{u_{l-1},u_l},
\end{align*}
where $Q=(Q_{z_1,z_2})_{z_1,z_2\in L_n}$ is the matrix given by
\[
Q_{z_1,z_2}=q(z_1,z_2)\ONE_{[z_1],[z_2]\subset [v]}.
\]
Therefore, the l.h.s.\ of  \eqref{eq:checking_Markov_property_1} equals
\[
P^{v}(t_1,\pi_nx,u_1)P^{v}(t_2-t_1,u_1,u_2)\ldots P^{v}(t_l-t_{l-1},u_{l-1},u_l),
\]
where
\begin{equation}
\label{eq:confined_transition_probabilities}
P^{v}(s,z_1,z_2)=\exp(q(v,v)s)\exp\{sQ\}_{z_1,z_2},\quad z_1,z_2\in L_n, [z_1],[z_2]\subset [v].
\end{equation}
Clearly, $P^{v}(s,z_1,z_2)$ is a Markov transition matrix, and it does not depend on $T$ which completes
the proof of the lemma.\epf 

The lemma above means that our Markov process conditioned on the fact that
it stays within $[v]$ up to time $T$ is also a Markov process  with transition probabilities
that do not depend on $T$. Therefore, we can consistently define this process up to infinite time. We denote the
resulting measure on infinite paths in $[v]$ by  $P^{v}_x$. The collection of these measures for all $x\in C$ is
a Markov family.

Let us now give a precise notion of self-similarity. 
We say that the Markov family $(P_x)_{x\in C}$ is self-similar, if for any $n$ there is a number $\alpha_n$ such that for
every $v\in L_n$ and every $x\in [v]$, and any map $h$ realizing
the similarity between $[v]$ and~$C$,
\begin{equation}
P_x^v\{X_{t_1}\in B_1,\ldots, X_{t_l}\in B_l\}=P_{h(x)}\{X_{\alpha_n t_1}\in h(B_1),\ldots, X_{\alpha_n t_l}\in h(B_l)\}.
\end{equation}

We shall say that a Markov family (as well as the associated Markov process, transition function, and semigroup) on $C$ is SSS (symmetric and self-similar) if it is stochastically continuous, Feller, isometry invariant, self-similar, and the trajectories of the associated Markov process are a.s.-c\`adl\`ag.

Suppose that the Markov family $P_x$ is SSS.
Let $h$ be a similarity map between~$[v]$ and~$C$. Then under $P^v_x$, $h(X)$ is a Markov process
that inherits the stochastic continuity, Feller property and isometry invariance from the original Markov family.
Therefore, all the above reasoning
for the Markov family $(P_x)$ applies to $h(X)$. In particular
the distribution of $h(X)$  under $(P_x^v)$ is completely determined by rates  $q^{v}_n$ that are defined for $h(X)$ under $(P_x^v)$ in the same way as
the rates $q_n$ are defined for $X$ under $(P_x)$ (these rates do not depend on $T$ as well).

Due to \eqref{eq:confined_transition_probabilities}, the jump rates for the confined process~$X$ under~$(P^v_x)$ are
\begin{equation}
\label{eq:jump_rates_for_conditioned_process}
q^v(z_1,z_2)=\lim_{t\downarrow 0}\frac{P^{v}(t,z_1,z_2)-\delta_{z_1,z_2}}{t}=q(z_1,z_2)+\delta_{z_1,z_2}q(v,v).
\end{equation}
This means that the process confined to $[v]$ can be viewed as the original process except the jumps out of $[v]$
are prohibited or ignored. 
Taking $v=(0)$ and the left $1$-shift on sequences $(0,x_2,x_3\ldots)$ for $h$, we get for the jump rates of $h(X)$ under~$(P^v_x)$:
\begin{equation}
q^{(0)}_n=q_{n+1},\quad n\in\N.
\label{eq:rough_jump_rates_for_conditioned_process} 
\end{equation}
By the self-similarity hypothesis we must have
\begin{equation}
q^{(0)}_n=\alpha_1 q_n,\quad n\in\N.
\label{eq:compare_new_jump_rates_to_the_old} 
\end{equation}
Comparing \eqref{eq:rough_jump_rates_for_conditioned_process} and \eqref{eq:compare_new_jump_rates_to_the_old}, we see that
\[
q_{n+1}=\alpha_1 q_n,,\quad n\in\N, 
\]
so that
\begin{equation}
q_n=\alpha_1^{n-1} q_1,\quad n\in\N.
\label{eq:geometric_progression}
\end{equation}

\begin{theorem}
Suppose a Markov family on $(C,d)$ is SSS. Then there are numbers $\theta,\gamma\ge0$ such that for every cylindric function $f=h\circ \pi_n$, the generator $Af$ is given by
\begin{equation}
Af(x)=\gamma\sum_{k=1}^n \theta^{k}(\langle h\rangle_{n,k,x}-h(\pi_n x)),
\label{eq:pregenerator}
\end{equation}
where
\[
 \langle h\rangle_{n,k,x}= 2^{-(n-k)}\sum_{v\in L_n:c(x,v)=k} h(v).
\]
\end{theorem}
\bpf Set $\theta=\alpha_1$, $\gamma=q_1/\alpha_1$ ($\gamma=0$ if $\alpha_1=0$), and use \eqref{eq:generator1}, \eqref{eq:jump_rates_uniformly_spread}, and
\eqref{eq:geometric_progression}.\epf

Let us recall that a linear operator $A$ defined on a vector subspace $\Dc$ of $\C(C)$ is called
a Markov pregenerator (see \cite[Chapter~1, Definition 2.1]{Liggett:MR776231}), if
\begin{enumerate}
 \item $1\in\Dc$ and $A1=0$.
 \item $\Dc$ is dense in $\C(C)$.
 \item \label{it:dissipativity} If $f\in\Dc$, $\mu\ge 0$ and $f-\mu Af =g$, then
\[
 \min_{x\in C} f(x)\ge\min_{x\in C} g(x).
\]
\end{enumerate}

\begin{lemma}\label{lm:pregenerator} The operator $A$ defined on cylindric functions via \eqref{eq:pregenerator} is
a Markov pregenerator. 
\end{lemma}
\bpf First two properties are obvious, and the third one follows 
from~\cite[Chapter 1, Proposition 2.2]{Liggett:MR776231} since
$A$ satisfies the following easily verifiable condition: if $f\in \Dc$ and
$f(x^*)=\min_{x\in C}f(x)$, then $Af(x^*)\ge 0$. \epf

Lemma~\ref{lm:pregenerator} implies that the closure of $A$ denoted by $\bar A$ is also a well-defined closed Markov pregenerator, see \cite[Chapter 1, Proposition 2.5]{Liggett:MR776231}.

\begin{lemma}
The operator $\bar A$ is a Markov generator, i.e., it is a closed Markov 
pregenerator satisfying 
\begin{equation}
\label{eq:image=everything}
\Rc(I-\mu A)=\C(C),\quad\mbox{for all}\ \mu>0.
\end{equation}
\end{lemma}
\bpf We need only to show~\eqref{eq:image=everything}. It is easy to see that
for any cylindric function $h$, there is a cylindric function $f$ such that
$f-\mu \bar A f=h$ (this can be derived from the fact that $\pi_n x$ is a Markov process for every $n$, and the one-to-one correspondence between Markov semigroups and Markov generators given by the Hille--Iosida
theorem). The lemma follows since the set of cylindric functions is dense, and $\Rc(I-\mu \bar A)$
is closed  (see \cite[Chapter 1, Proposition 2.6]{Liggett:MR776231}). \epf

We can now summarize the above.
\begin{theorem}\label{th:necessary_cond} If a Markov family $(P_x)_{x\in C}$ is SSS,
then the generator of the associated Markov semigroup coincides with the closure of the operator $A$ defined on cylindric functions 
via~\eqref{eq:pregenerator} for some $\gamma,\theta\ge0$.
\end{theorem}\label{th:sufficient_cond}
\bpf The result follows now from the Hille-Iosida theorem (see e.g.~\cite[Chapter 1, Theorem 2.9]{Liggett:MR776231}) which establishes a one-to-one
correspondence between Markov generators and Markov semigroups.\epf

Theorem~\ref{th:necessary_cond} gives a necessary condition for a 
Markov semigroup to be SSS. The next result shows that this condition is, in fact, sufficent. 

\begin{theorem} For any $\gamma,\theta\ge 0$ there is a unique Markov family with Markov generator coinciding with $A$ defined on cylindric functions  via~\eqref{eq:pregenerator}. That Markov family is SSS (with
scaling parameter given by $\alpha_n=\theta^{n}$, $n\in\N$). 
\end{theorem}
\bpf The existence-uniqueness and the Feller property follows from the Hille--Iosida theorem.
The isometry invariance follows from that of $A$. The stochastic continuity follows from
\[
P(t,x,[v])=1-t\gamma\sum_{k=1}^n\theta^{k}+o(t),\quad t\to 0,
\]
for any $n\in\N$, any $v\in L_n$, and any $x\in[v]$.
In particular, a c\`adl\`ag version of the canonical process exists.

For the self-similarity, we must take
any $n\in\N$, $v\in L_n$, $x\in [v]$, and consider the process $X$ emitted from $x$ under the condition that
it stays within $[v]$. Due to Lemma~\ref{lm:restriction_is_markov}, the conditioned process is Markov, and so is $h(X)$ under~$P^v$,
where $h$ is a similarity map between $[v]$ and $C$. Computing the transition probabilities and jump rates for this process:
\[
 q^v_{k}=q_{n+k}=\gamma\theta^{n+k}=\theta^n\gamma\theta^{k}=\theta^n q_k,
\]
so that $h(X)$ under $P$ has the same distribution as $X$ under $P^v$ and time change $t\to\theta^nt$. The
proof is complete. \epf

Theorems~\ref{th:necessary_cond} and \ref{th:sufficient_cond} give a complete characterization of SSS processes, the analogues of the Wiener
process on $(C,d)$ via their Markov generators. SSS processes are naturally parametrized by $\gamma$ and $\theta$. We shall write
$SSS(\gamma,\theta)$ to denote the SSS process with parameters $\gamma,\theta$. It is clear that $\gamma$ is responsible for uniform time changes, and it is often sufficient to study the case $\gamma=1$, since by a simple time rescaling
one can obtain the process with any given $\gamma$. However, the self-similarity parameter $\theta$ is essential, and there are
qualitative differences between processes with different values of $\theta$. 

It is important to stress that the class of SSS processes that we describe is the same as
the class of processes obtained in~\cite{Pearson-Bellisard} via the noncommutative geometry approach. The correspondence
can be established by noticing that
our scaling factor $\theta$ is equal to $3^{s_0+2-s}$ in the notation of~\cite{Pearson-Bellisard}, where $s_0 =\frac{\ln 2}{\ln 3}$ (the box dimension of $C$).

Since SSS processes on the Cantor set play the role of symmetric diffusion, the generators $A$ play the role
of the Laplacian. We shall study spectral properties of $A$ in Section~\ref{sec:spectral_structure}.

\section{Explicit construction and jump statistics}\label{sec:explicit_construction}

We begin with an explicit construction of  $SSS(\gamma,\theta)$.
Let $(\xi_{kj})_{k,j=1}^\infty$ be a family of independent random variables such that for each $k\in\N$,
$(\xi_{k1},\xi_{k2},\ldots)$ are exponentially distributed with parameter~$\gamma\theta^{k}$.
Let $S_{kn}=\xi_{k1}+\ldots+\xi_{kn}$. Clearly, $N_k(t)=\max\{n: S_{kn}\le t\}$, $t\ge0$,
is a c\`adl\`ag Poisson process with intensity $\gamma\theta^{k-1}$. We shall say that there is a jump at level $k$ at time $\tau$ if $N_k(\tau)=N_k(\tau-)+1$.
Let $\tilde N_k(t)=N_1(t)+\ldots+N_{k-1}$ be the Poisson process that counts the jumps of Poisson processes at all levels below $k$, i.e. at levels $1,2,\ldots,k-1$. 

To define our Markov process $X$ we shall also need a family of i.i.d. $\frac{1}{2}$-Bernoulli random variables $(\kap_{kj})_{k,j=1}^\infty$ independent
of the Poisson processes described above.

We set $X(0)=x=(x_1,x_2,\ldots)$ and let the evolution of the $k$-th coordinate $X_k$ to be defined by the following rules:
\begin{enumerate} 
\item $X_k$ stays constant
while the processes $N_k$ and $\tilde N_{k}$ are constant. 
\item If at time $\tau$ the process $N_k$ makes a jump, $X_k$ also makes a jump so that
$X_k(\tau)=1-X_k(\tau-)$. 
\item If at time $\tau$ the process $\tilde N_k$ makes a jump, then $X_k(\tau)=\varkappa_{k \tilde N_k(\tau)}$ no matter what the value
of $X_k(\tau-)$ was.
\end{enumerate}

In other words, when a jump occurs at a level $k$, $X_1,\ldots,X_{k-1}$ do not change,  $X_k$ gets flipped, and $(X_{k+1},X_{k+2},\ldots)$ are
re-initialized according to the $\frac{1}{2}$-Bernoulli product measure. We exclude the event of probability $0$ on which two jumps happen at the same time. The process $(X_1,\ldots,X_n)$ makes finitely many jumps in finite time for any finite $n$.

It is easy to see that this procedure uniquely defines a Markov process with pregenerator described in the last section.

We see that the value $\theta=1$ is critical. If $\theta<1$ then with probability $1$, $X$
makes finitely many jumps in a finite time. However, if $\theta\ge 1$, then with probability $1$, $X$ makes
infinitely many jumps in a finite time.

The value $\theta=1$ is also special due to the following: if $\theta=1$ then for any $n\in\N$ and any  $v\in L_n$, the distribution of the restricted process $h(X)$ under $P^v$ coincides precisely
with that of $X$ (no time change is needed). Notice also that in terms of this model our conditioning means that the process $\tilde N_n+N_n$ makes no jumps.

\section{Spectral structure of the generator}
\label{sec:spectral_structure}
In this section we fix the values of $\theta>0$ and $\gamma=1$ and study eigenvalues and eigenvectors of the infinitesimal operator $A$. The eigenvectors will be given by Haar function that we proceed to introduce. For any $n$ and any $v\in L_n$, we define
\[
\psi_v=2^{n/2}(\chi_{v0}-\chi_{v1}),
\]
where $\chi_{u}=\chi_{[u]}$ denotes the characteristic function (or indicator) of $[u]$ for any $u$. Notice that $\psi$ with no indices
denotes $\chi_0-\chi_1$.
\begin{theorem}\label{th:spectral}
\begin{enumerate}
\item
Eigenvalues of $A$ are given by $\lambda_0=0$, and
\begin{align*}
\lambda_n&=-\sum_{k=1}^{n-1}\theta^k-2\theta^n,\quad n\in\N\\
         &\left(=-\frac{2\theta^{n+1}-\theta^n-\theta}{\theta-1}\ \text{if}\ \theta\ne 1\right).
\end{align*}
\item The unique (up to a multiplicative constant) eigenfunction associated to $\lambda_0$ is $1$. For $n\in\N$, $\lambda_n$ has multiplicity $2^{n-1}$, and its eigenspace
is spanned by $M_n=\{\psi_v:\ v\in L_{n-1}\}$.
\item The eigenfunctions described above form a complete basis in $\C(C)$ and~$L^2(C)$. This basis is orthonormal in $L^2(C,b)$, where
$b$ is the $1/2$-Bernoulli measure on $C$ (also known as the Cantor measure).
\item If $\theta<1$, $A$ is bounded with spectral radius equal to
\[
r(A)=\sum_{k=1}^{\infty}\theta^k = \frac{\theta}{1-\theta}.
\]
Moreover, $A+r(A)I$ is compact where $I$ is the identity operator, so that $A$ is a compact perturbation of a multiple of the identity.
\item If $\theta\ge 1$ then $A$ has compact resolvent.
\end{enumerate}
\end{theorem}
This theorem is a result of straightforward computations and arguments. We omit the proof.

\section{Unique ergodicity, exponential mixing}\label{sec:ergodicity}
\begin{theorem} \begin{enumerate}\item For any $\gamma,\theta>0$, the $\frac{1}{2}$-Bernoulli product measure $b=\mu^\N$ on $C$ is a unique invariant measure for $SSS(\gamma,\mu)$. \item There are constants $K,\nu>0$, such that for any point
$x\in C$,
\[
|P(t,x,\cdot)-b(\cdot)|_{TV}\le K e^{-\nu t},
\]
 where $|\cdot|_{TV}$ denotes the total variation norm.
\end{enumerate}
\end{theorem}
\bpf The Bernoulli measure $b$ is clearly invariant for SSS since it is invariant under isometries.
The uniqueness of the invariant measure follows from the second statement of the theorem. The latter can be proved by exploiting the spectral gap provided by Theorem~\ref{th:spectral}, but we choose another simple method instead.

Let us work with the explicit model introduced in Section~\ref{sec:explicit_construction}.
It is easy to see that for every $N\ge1$, the distribution of $X(t)$ conditioned on $\{N_1(t)=N\}$
is given by $\delta_{(N+x_1)\imod 2}\times b$. This implies that for any set $B$,
\begin{align*}
P(t,x,B)=&\sum_{N=1}^\infty \Pp\{N_1(t)=N\}(\delta_{(N+x_1)\imod 2}\times b)(B)+\beta(t,B)\\
        =&(\delta_{(x_1+1)\imod2}\times b) (B)e^{-\gamma\theta t}\sum_{m=1}^{\infty}\frac{(\gamma\theta t)^{2m-1}}{(2m-1)!}\\
	&+(\delta_{x_1\imod2}\times b) (B)e^{-\gamma\theta t}\sum_{m=1}^{\infty}\frac{(\gamma\theta t)^{2m}}{(2m)!}+\beta(t,B)\\
	=&(\delta_{(x_1+1)\imod2}\times b)(B)e^{-\gamma\theta t}\sinh(\gamma\theta t) \\
	&+(\delta_{x_1\imod2}\times b)(B)e^{-\gamma\theta t}(\cosh(\gamma\theta t)-1)+\beta(t,B)\\
	=&\frac{1}{2}(\delta_{(x_1+1)\imod2}\times b)(B)(1-e^{-2\gamma\theta t}) \\
	&+\frac{1}{2}(\delta_{x_1\imod2}\times b)(B)(1+e^{-2\gamma\theta t}-2e^{-\gamma\theta t})+\beta(t,B),\\
	=&b(B)-\frac{1}{2}(\delta_{(x_1+1)\imod2}\times b)(B)e^{-2\gamma\theta t}\\
	&+\frac{1}{2}(\delta_{x_1\imod2}\times b)(B)(e^{-2\gamma\theta t}-2e^{-\gamma\theta t})+\beta(t,B),\\
\end{align*}
where 
\[
0\le\beta(t,B)\le \Pp\{N_1(t)=0\}=e^{-\gamma\theta t}.
\]
So the second statement of the theorem follows with $K=\frac{3}{2}$, $\nu=\gamma\theta$.
\epf

\section{Displacement moments}\label{sec:moments}
In this section we study the behavior of the moments of displacement of the process $X$ given by
\[
M_r(t)=\E_x d^r(x,X(t))=\sum_{k=1}^\infty 3^{-rk}P_x(c(x,X(t))=k),
\]
as $t\to 0$.
Due to the isometry invariance it is sufficient to consider $x=\bar 0=(0,0,\ldots)$.
\begin{theorem}\label{thm:moments} Let $r>0$.
 \begin{enumerate}
  \item If $3^r>\theta$, then
\[
M_r(t)=\sum_{k=1}^\infty 3^{-rk}\theta^k(t+o(t)),\quad t\downarrow 0.
\]
 \item If $3^r<\theta$, then
\[0<\liminf_{t\downarrow 0}\frac{M_r(t)}{t^{\frac{\ln 3}{\ln\theta}r}}\le  
\limsup_{t\downarrow 0}\frac{M_r(t)}{t^{\frac{\ln 3}{\ln\theta}r}}<\infty.\]
 \end{enumerate}
\end{theorem}
\begin{remark} It is shown in \cite{Pearson-Bellisard} that if
$3^r=\theta$, then $M_r(t)$ behaves as $t\ln(1/t)$ as $t\to0$.
\end{remark}


\bpf[Proof of Theorem~\ref{thm:moments}]
For any $k\in\N$, 
\begin{align*}
P(x,t,[v_k(x)])&=\int_{y\ne x} \chi_{[v_k(x)]}(y) P(x,t,dy)\\
               &=S^t \chi_{[v_k(x)]} (x)\\
	       &=\langle \chi_{[v_k(x)]},1\rangle  \chi_{[v_k(x)]}(x) + \sum_{i=0}^{k-1} e^{\lambda_{i+1} t}\sum_{u\in L_i} \langle \chi_{[v_k(x)]},\psi_u\rangle \psi_u(x),
\end{align*}
where $\langle\cdot,\cdot\rangle$ denotes the inner product in the Hilbert space $L^2(C,b)$.
The first term in the r.h.s.\ is $0$, and for each $i$ in the sum above, the only nonzero contribution in the sum comes from
\[
\langle \chi_{[v_k(x)]},\psi_{\pi_{i}x}\rangle=
\begin{cases}
2^{\frac{i}{2}-k},& i<k-1\\
-2^{-\frac{k+1}{2}}&     i=k-1,
\end{cases}
\]
so that
\[
P(x,t,[v_k(x)])=\sum_{i=0}^{k-2} e^{\lambda_{i+1} t} 2^{\frac{i}{2}-k}\bar \psi_i(\bar 0) -e^{\lambda_{k} t} 
2^{-\frac{k+1}{2}}\bar \psi_{k-1}(\bar 0),
\]
where
\[
\bar \psi_i=\psi_{\underbrace{\scriptstyle 00\ldots0}_{i}}.
\]
Since $\psi_i(\bar 0)=2^{i/2}$, we obtain
\[
P(x,t,[v_k(x)])=\sum_{i=0}^{k-2} 2^{i-k}e^{\lambda_{i+1} t}  -2^{-1}e^{\lambda_{k} t},
\]
and
\[
M_r(t)=\sum_{k=1}^{\infty} 3^{-kr}\left(\sum_{i=0}^{k-2} 2^{i-k} e^{\lambda_{i+1} t} -2^{-1}e^{\lambda_{k} t}\right)
\]
If $3^r>\theta$, then the series
\[
\sum_{k=1}^{\infty} 3^{-kr}\left(\sum_{i=0}^{k-2} \lambda_{i+1} 2^{i-k} e^{\lambda_{i+1} t} -2^{-1}\lambda_{k}
e^{\lambda_{k} t}\right) 
\]
converges uniformly in $t\ge 0$. Therefore, it represents $\partial_t M_r(t)$. Evaluating this series at $0$ produces
\[
\frac{d^+}{dt} M_r(t)=\sum_{k=1}^{\infty} 3^{-kr}\left(\sum_{i=0}^{k-2} \lambda_{i+1} 2^{i-k} -2^{-1}\lambda_{k}\right).
\]
One can show by induction that
\[
\sum_{i=0}^{k-2} \lambda_{i+1} 2^{i-k} -2^{-1}\lambda_{k}=\theta^k,\quad k\in\N, 
\]
so that in this case
\[
\frac{d^+}{dt}M_r(t)=\sum_{k=1}^\infty 3^{-rk}\theta^k,
\]
and the first statement of the theorem follows.

Let now $\theta>3^r>1$. Consider a sequence
$t_n=R_n\theta^{-n}$, $n\in\N$ with $R_n\in[1,\theta]$.

We have
\begin{equation}
\label{eq:supercritical_moments_aux1}
 M_r(t_n)=\sum_{k=1}^{\infty}a_{k,n}=\sum_{m=-n+1}^{\infty}a_{n+m,n},
\end{equation}
where
\[
a_{k,n}=3^{-kr}\left(\sum_{i=0}^{k-2} 2^{i-k} e^{\lambda_{i+1} t_n} -2^{-1}e^{\lambda_{k} t_n}\right)>0.
\]
Taking $k=m+n$, we obtain
\begin{align*}
\frac{a_{m+n,n}}{t^{\frac{\ln 3}{\ln\theta}r}}&=\frac{3^{-(n+m)r}}{R_n^{\frac{\ln 3}{\ln\theta}r}3^{-nr} }\left(\sum_{i=0}^{n+m-2} 2^{i-(n+m)} e^{\lambda_{i+1} R_n\theta^{-n}} 
-2^{-1}e^{ \lambda_{n+m} R_n\theta^{-n}}\right)\\
&=R_n^{-\frac{\ln 3}{\ln\theta}r} 3^{-mr}\left(\sum_{l=1}^{n+m-1} 2^{-1-l} 
e^{\lambda_{n+m-l}R_n\theta^{-n}} 
-2^{-1}e^{\lambda_{n+m}R_n\theta^{-n} }\right),
\end{align*}
where we used a change of variables $l=n+m-1-i$. Therefore, for $m\ge 0$,
\begin{equation}
\label{eq:a_m_n_3}
\frac{a_{m+n,n}}{t^{\frac{\ln 3}{\ln\theta}r}}\le 3^{-mr} \left(\sum_{l=1}^{n+m-1} 2^{-1-l}+2^{-1}\right)\le 3^{-mr}.  
\end{equation}
If $m\le 0$, then we use another estimate:
\begin{equation}\label{eq:a_m_n}
\frac{a_{m+n,n}}{t^{\frac{\ln 3}{\ln\theta}r}}
\le 3^{-mr}\left(1-e^{\lambda_{n+m}R_n\theta^{-n} }\right).
\end{equation}
Since
\[
 |\lambda_{n+m}R_n\theta^{-n}|\le \frac{2\theta^{n+m+1}-\theta^{n+m}-\theta}{\theta-1}\theta^{-n+1}\le K_1\theta^m,
\]
for some $K_1>0$, all $n>0$, $-n+1\ge m\le0$, inequality \eqref{eq:a_m_n} can be continued as
\begin{equation}\label{eq:a_m_n_2}
 \frac{a_{m+n,n}}{t^{\frac{\ln 3}{\ln\theta}r}}
\le K_2 3^{-mr}\theta^m,
\end{equation}
for some $K_2>0$ and all $m\le0$. The upper estimate in the theorem follows now 
from~\eqref{eq:a_m_n_3} and~\eqref{eq:a_m_n_2}. To prove the lower bound, take $m=0$:

\begin{align*}
 \frac{M_r(t)}{t^{\frac{\ln 3}{\ln\theta}r}}\ge\frac{a_{n,n}}{t^{\frac{\ln 3}{\ln\theta}r}}
\ge R_n^{-\frac{\ln 3}{\ln\theta}r} \left(\sum_{l=1}^{n-1} 2^{-1-l} 
e^{\lambda_{n-l}R_n\theta^{-n}} 
-2^{-1}e^{\lambda_{n}R_n\theta^{-n} }\right),
\end{align*}
Notice that
\begin{equation}\label{eq:conv_to_nu}
\lambda_{n-l}\theta^{-n}= -\frac{2\theta^{-l+1}-\theta^{-l}-\theta^{1-n}}{\theta-1}\ \to\ -\nu\theta^{-l},\quad n\to\infty,
\end{equation}
where $\nu=\frac{2\theta-1}{\theta-1}$. So,
\begin{align*}
 \frac{M_r(t)}{t^{\frac{\ln 3}{\ln\theta}r}}
&\ge \theta^{-\frac{\ln 3}{\ln\theta}r} \left(\sum_{l=1}^{n-1} 2^{-1-l} 
e^{-R_n\nu\theta^{-l}}
-2^{-1} e^{-R_n\nu}\right)\\
&+\theta^{-\frac{\ln 3}{\ln\theta}r} \left(\sum_{l=1}^{n-1} 2^{-1-l} 
(e^{\lambda_{n-l}R_n\theta^{-n}}-e^{-R_n\nu\theta^{-l}}) 
-2^{-1}(e^{\lambda_{n}R_n\theta^{-n}}- e^{-R_n\nu})\right).
\end{align*}
Due to~\eqref{eq:conv_to_nu}, the second sum converges to zero, and we have
\[
 \liminf_{n\to\infty}\frac{M_r(t)}{t^{\frac{\ln 3}{\ln\theta}r}}\ge\theta^{-\frac{\ln 3}{\ln\theta}r} \inf_{R\in[1,\theta]}
 \left(\sum_{l=1}^{\infty} 2^{-1-l} 
e^{-R\nu\theta^{-l}}
-2^{-1} e^{-R\nu}\right)>0,
\]
and the proof is completed.
\epf

\bibliographystyle{alpha}
\bibliography{cantor}
\end{document}